\documentclass[12pt]{amsart}
\usepackage{amsfonts}
\usepackage{amsmath,amssymb,amsthm}
\usepackage{amstext}
\usepackage{verbatim}
\input{epsf}
\newtheorem{thm}{Theorem}[section]

\newtheorem{prob}[thm]{Problem}

\theoremstyle{definition}

\theoremstyle{remark}

\begin{document}

\title[Algorithm for composition of inverse problems with exclusive disjunction]
{An algorithm for composition of inverse problems with exclusive disjunction as a
logical structure in the conclusion}
\author{Julia Ninova and  Vesselka Mihova}%
\address{Faculty of Mathematics and Informatics University of Sofia
e-mail: julianinova@hotmail.com}%
\address{Faculty of Mathematics and Informatics University of Sofia
e-mail: mihova@fmi.uni-sofia.bg}


\maketitle \thispagestyle{empty}

\begin{abstract}
By suitable examples we illustrate an algorithm for composition of inverse problems.
\end{abstract}

\section{introduction}
Let $\,t, p_1, p_2, r$ be given statements. We deal with a generalization of the formal logical rule \cite{S}
$$(p_1\rightarrow r)\,\wedge\,(p_2\rightarrow r)\;\Leftrightarrow\; p_1\vee p_2\,\rightarrow\, r.$$

In \cite{MN1} and \cite{MN2} we prove, clarify and use the equivalence
$$(t\wedge p_1\rightarrow r)\,\wedge\,(t\wedge p_2\rightarrow r)\;\Leftrightarrow\;
t\wedge (p_1\vee p_2)\,\rightarrow\, r. \leqno (*)$$

It gives an algorithm for composition of inverse problems with a given
logical structure that is based on the steps below.
\begin{itemize}
\item[-] Formulating  and proving  \emph{generating problems} with logical structures
of the statements as those at the left hand side of (*).

\item[-] Formulating a problem with a logical structure of the statement
$\;t\wedge (p_1\vee p_2)\rightarrow r.$

\item[-] Formulating an \emph{inverse problem}
with a logical structure $\,t\wedge r \rightarrow  p_1\vee p_2$.
\end{itemize}

A technology for composition of equivalent and inverse
 problems  is appropriate for training of mathematics students and teachers.
\vskip 2mm

In what follows we illustrate the above algorithm by two groups of problems.

\section{Group of problems I}

The logical statements used for the formulation of the problems in this group are

\begin{itemize}
\item[] $t:\; \{$ \emph{In the plane are given a $\triangle\, ABC$, a straight line $g\ni C$ and an arbitrary
point} $M\in g, M\neq C.\}$
\vskip 2mm

\item[] $p_1:\;\{$ \emph{The straight line $g$ is a median in} $\triangle\, ABC .\}$
\vskip 2mm

\item[] $p_2:\;\{ g\,\parallel\, AB\}$
\vskip 2mm

\item[] $r:\; \{ S_{\triangle\, AMC} = S_{\triangle\, BMC} \}$
\end{itemize}

\begin{figure}[h t b]
\epsfxsize=8cm \centerline{\epsfbox{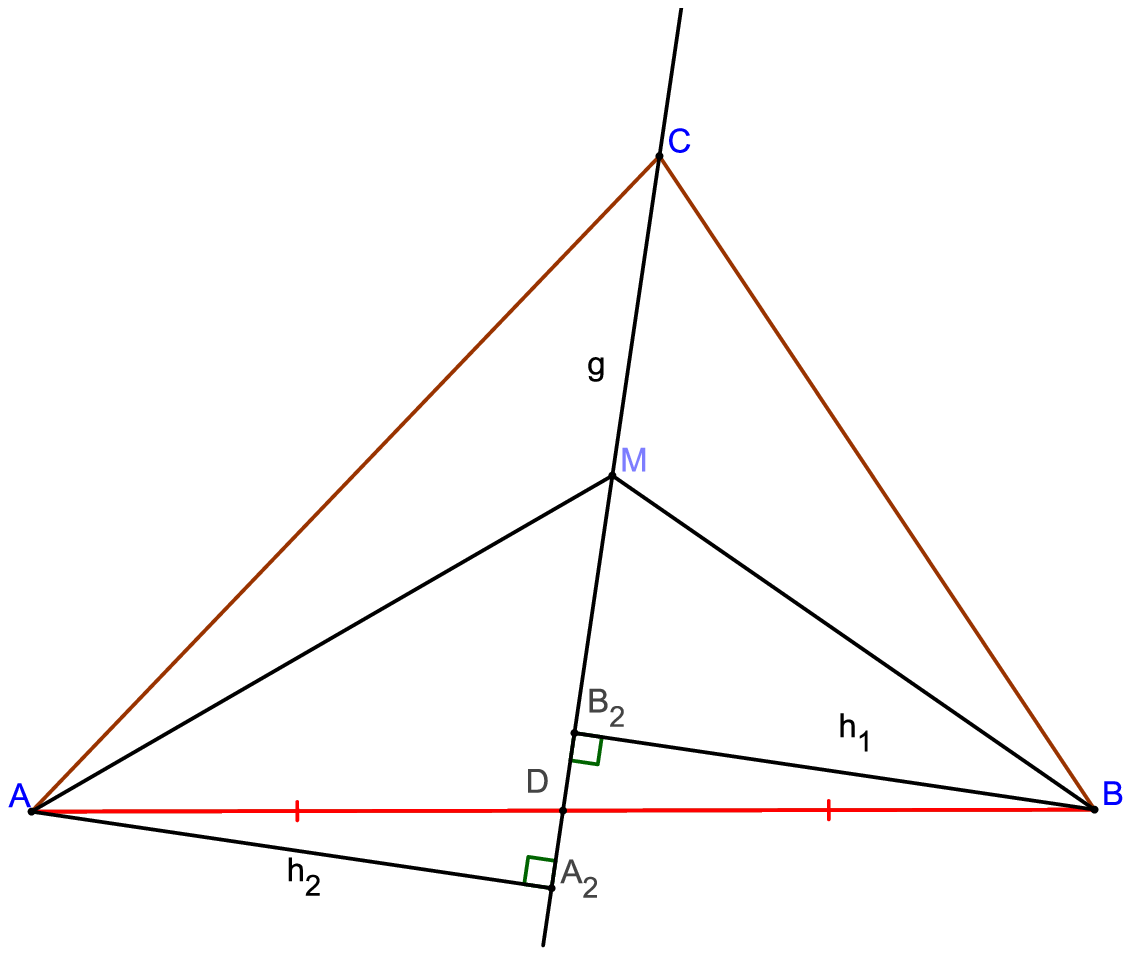}}
\end{figure}

\begin{prob}
Let in a $\triangle\, ABC$ the straight line $g$ passes through
$C$ and is a median. Prove that for any point $M$ on $g$ $(M\neq
C)$ the triangles $\triangle\, AMC$ and $\triangle\, BMC$ are
equal in area.
\end{prob}

This problem has a logical structure $t\wedge p_1 \rightarrow r$.
\vskip 1mm

\emph{Proof.} Since the straight line $g$ passes through $C$ and is a median in a $\triangle\, ABC$,
then $AB \cap g\neq \emptyset$.
Let $AB \cap g= D$ (fig. 1). Hence $AD = BD$.

Let $h_1$ and $h_2$ be the altitudes in $\triangle\, BDC$ and $\triangle\, ADC$
through the vertices $B$ and $A$ respectively. Then
$$\frac{S_{\triangle\, AMC}}{S_{\triangle\, BMC}}=\frac{CM.h_2}{CM.h_1}=\frac{h_2}{h_1}=
\frac{MD.h_2}{MD.h_1}=\frac{S_{\triangle\, AMD}}{S_{\triangle\, BMD}}=1.$$
\hfill{$\square$}
\vskip 2mm

\begin{prob}
Let in $\triangle\, ABC$ the straight line $g$ passes
through $C$ and is parallel to the straight line $AB$. Prove that
for any point $M$ on $g$ $(M\neq C)$ the triangles $\triangle\,
AMC$ and $\triangle\, BMC$ are equal in area.
\end{prob}

This problem has a logical structure $t\wedge p_2 \rightarrow r$.
\vskip 1mm

The proof follows immediately from Figure 2.

\begin{figure}[h t b]
\epsfxsize=8cm \centerline{\epsfbox{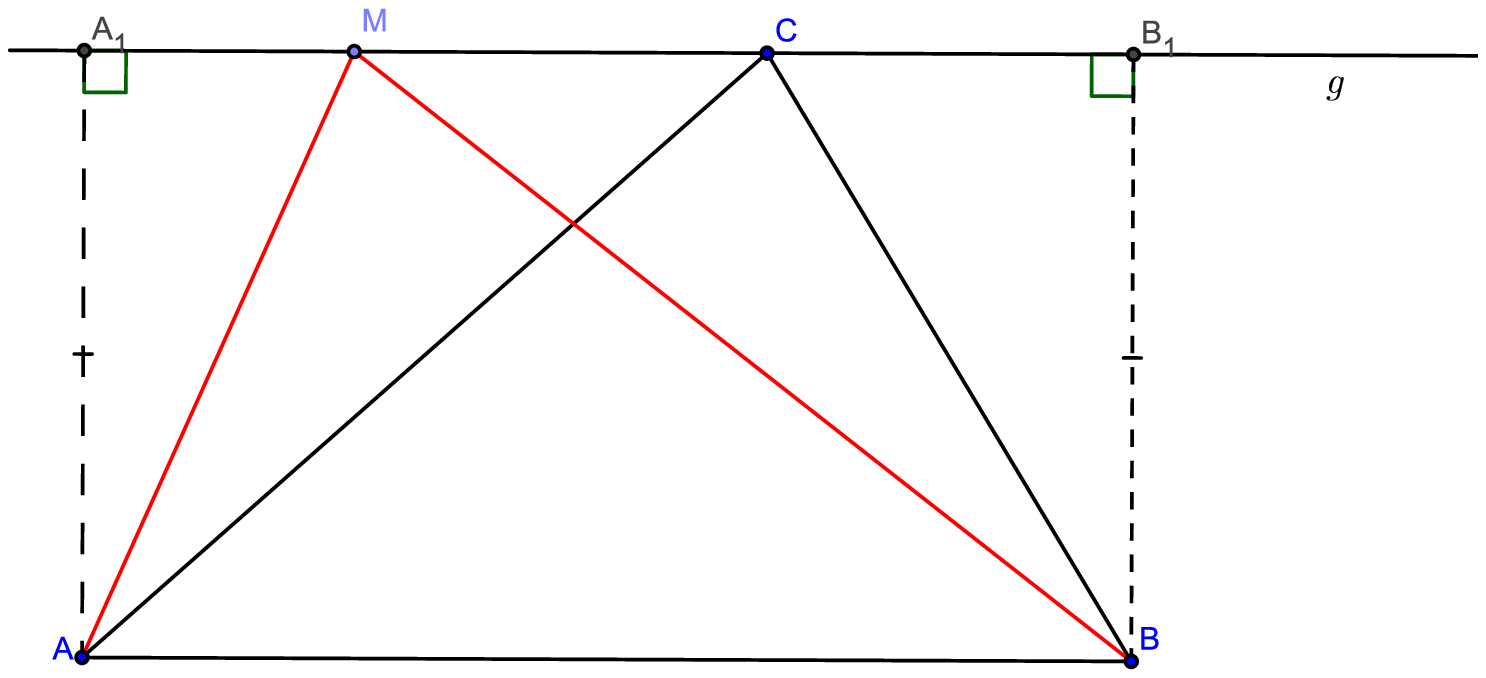}}
\end{figure}

According to (*) and the \emph{generating problems} 2.1 and 2.2, we can formulate a
problem with a logical structure $\;t\wedge (p_1\vee p_2)\rightarrow r.$

Since the statements $p_1$ and $p_2$ are mutually exclusive,
then $\;p_1\vee p_2\,\Leftrightarrow\, p_1\,\veebar\, p_2$. Under this condition problems
with structures $\;t\wedge (p_1\vee p_2)\rightarrow r$ and $\;t\wedge (p_1\veebar p_2)\rightarrow r$
are equivalent.
\vskip 1mm

We formulate and solve the following inverse problem with a logical structure
$$t\wedge r \rightarrow p_1\veebar p_2.$$

\begin{prob}
Let the point $M$ in the plane of $\triangle\, ABC$ is such that
$\triangle\, AMC$ and $\triangle\, BMC$ exist and are equal in area. Prove that the straight line $CM$ either
cuts the side $AB$ of $\triangle\, ABC$ at its middle point or is parallel to $AB$.
\end{prob}

\emph{Proof.} The triangles $\triangle\, AMC$ and $\triangle\, BMC$ have a common side $CM$ and are equal in area.
Hence the altitudes to their common side are equal, i. e. the vertexes $A$ and $B$ are equidistant from $CM$.

There are two possibilities for the location of the points $A$ and $B$ with respect to the straight line $CM$.
\begin{itemize}
\item[-] If $A$ and $B$ are situated at opposite sides of $CM$  (fig. 1), then
the point $D=AB \cap CM$  is the middle point of the side $AB$
($\triangle\, AA_2D\cong \triangle\, BB_2D$) and the straight line $CM$ is the median of  $\triangle\, ABC$ through $C$.

\item[-] If $A$ and $B$ are situated at one and the same side of $CM$ (fig. 2), then
they lie on a straight line parallel to $CM$, i. e. $CM \parallel AB$.
\end{itemize}
\hfill{$\square$}

\section{Group of problems II}
The logical statements used for the formulation of the problems in this group are
\begin{itemize}
\item[] $t:\; \{$ \emph{Let $ABCD$ be a quadrilateral with} $AC\cap BD=O.\}$
\vskip 2mm

\item[] $p_1:\;\{ AB \parallel CD\}$
\vskip 2mm

\item[] $p_2:\;\left\{ \displaystyle{\frac{AB}{CD}=1}\right\}$
\vskip 2mm

\item[] $r:\; \left\{ \displaystyle{\frac{AO}{OC}=\frac{BO}{OD}=\lambda} \right \}$
\end{itemize}

\begin{prob}
Let $ABCD$ be a trapezium with $AC\cap BD=O.$ Prove that
$$  \frac{AO}{OC}=\frac{BO}{OD}=\lambda,\, \lambda\neq 1. $$
\end{prob}

This problem has a logical structure $\,t\wedge (p_1\wedge \overline{p_2})\rightarrow r$.
\vskip 1mm

\begin{figure}[h t b]
\epsfxsize=8cm \centerline{\epsfbox{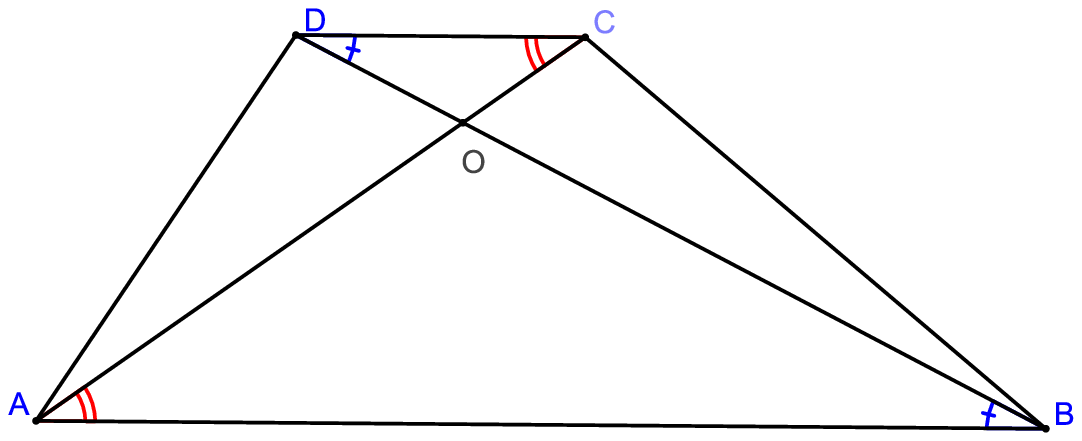}}
\end{figure}

\emph{Proof.} Since $\triangle \,AOB \sim \triangle\, COD$ (fig. 3) and $AB\neq CD$, then
$$\frac{AO}{OC}=\frac{BO}{OD}=\frac{AB}{CD}=\lambda,\,\lambda\neq 1. $$
\hfill{$\square$}

\begin{prob}
Let $ABCD$ be a parallelogram with $AC\cap BD=O.$ Prove that
$$  \frac{AO}{OC}=\frac{BO}{OD}=\lambda,\, \lambda= 1. $$
\end{prob}

This problem has a logical structure $\,t\wedge (p_1\wedge p_2)\rightarrow r$.
\vskip 1mm

\begin{figure}[h t b]
\epsfxsize=8cm \centerline{\epsfbox{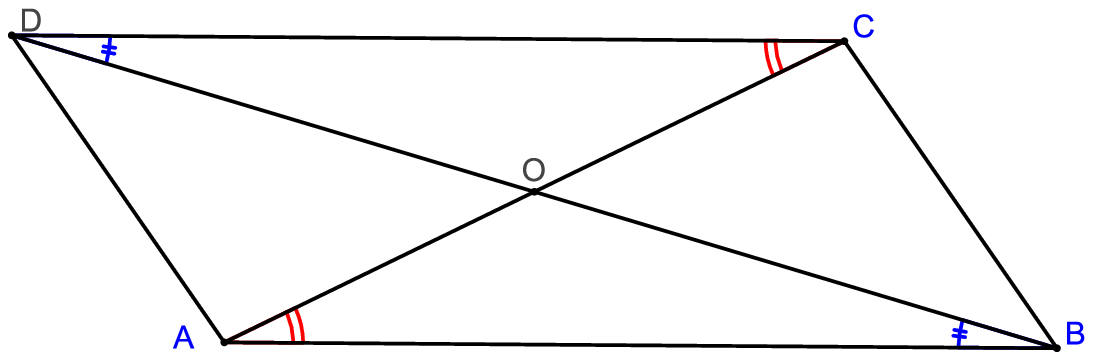}}
\end{figure}

\emph{Proof.} Since $AB=CD$, the proportion
$$\frac{AO}{OC}=\frac{BO}{OD}=\frac{AB}{CD}=\lambda,\,\lambda= 1 $$

\vskip 2mm
\noindent
follows from the congruence $\,\triangle \,AOB \cong \triangle\, COD$ (fig. 4).
\hfill{$\square$}
\vskip 2mm

According to (*) and the \emph{generating problems} 3.1 and 3.2, we can formulate a
problem with a logical structure $\;t\wedge ((p_1\wedge \overline{p_2})\vee (p_1\wedge p_2))\rightarrow r.$

Since the notion trapezium is not generic with respect to the notion parallelogram, then
problems with logical structures $\;t\wedge ((p_1\wedge \overline{p_2})\vee (p_1\wedge p_2))\rightarrow r$
and $\;t\wedge ((p_1\wedge \overline{p_2})\veebar (p_1\wedge p_2))\rightarrow r$ are equivalent.

We formulate and solve the following inverse problem with a logical structure
$$t\wedge r \rightarrow (p_1\wedge \overline{p_2})\veebar (p_1\wedge p_2).$$

\begin{prob}
Let $ABCD$ be a quadrilateral with $AC\cap BD=O.$ Prove that if
$$\frac{AO}{OC}=\frac{BO}{OD}=\lambda,$$
then the quadrilateral is either a trapezium or a parallelogram.
\end{prob}

\emph{Proof.} From the given proportionality and the equality $\, \angle \,AOB=\angle\, COD$
we conclude that $\,\triangle\,AOB\sim\triangle\,COD$.
\begin{figure}[h t b]
\epsfxsize=8cm \centerline{\epsfbox{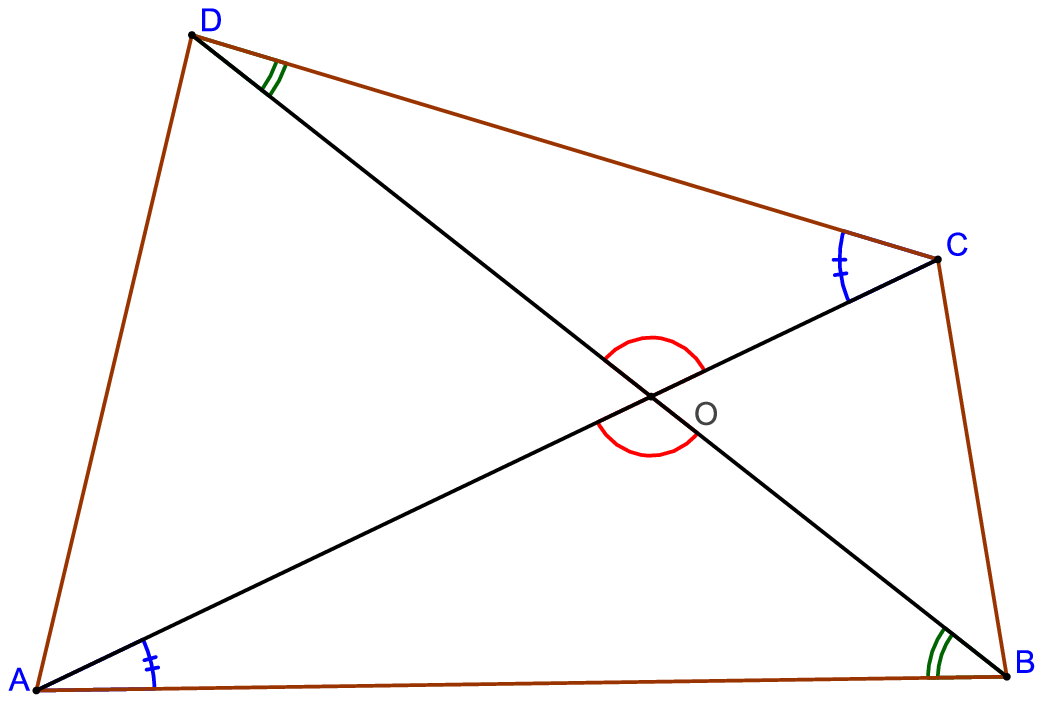}}
\end{figure}
Hence  $\,\displaystyle{\frac{AB}{CD}=\lambda}$,
$\, \angle \,OAB=\angle\, OCD$ and $\, \angle \,OBA=\angle\, ODC$, which implies
 the straight lines $AB$ and $CD$  are parallel.

There are two possibilities for the ratio coefficient $\lambda$: either $\lambda\neq 1$ or $\lambda= 1$.

\begin{itemize}
\item[-] If $\lambda\neq 1$  the quadrilateral $ABCD$ is a trapezium.
\vskip 1mm

\item[-] If $\lambda= 1$  the quadrilateral $ABCD$ is a parallelogram.
\end{itemize}
\hfill{$\square$}

\end{document}